\documentclass[11pt]{article}

\usepackage{amsmath}
\usepackage{amsfonts}
\usepackage{amssymb}
\usepackage{amsthm}
\usepackage{hyperref}
\usepackage{geometry}

\newtheorem{theorem}{Theorem}
\numberwithin{theorem}{section}
\newtheorem{proposition}[theorem]{Proposition}
\newtheorem{lemma}[theorem]{Lemma}
\newtheorem{remark}[theorem]{Remark}
\newtheorem{problem}[theorem]{Problem}
\newtheorem{corollary}[theorem]{Corollary}
\theoremstyle{definition}
\newtheorem{definition}[theorem]{Definition}
\newtheorem{example}[theorem]{Example}

\title{Algebraic Constraints for Linear Acyclic Causal Models}
\author{Cole Gigliotti and Elina Robeva}
\date{July 1, 2025}

\usepackage{graphicx} 
\usepackage{subcaption}

\renewcommand{\cal}{\mathcal}
\newcommand{\bb}{\mathbb}

\usepackage{tikz}
\usetikzlibrary{arrows}

\usepackage{natbib}
\usepackage{xcolor}
\usepackage{blkarray}
\usepackage{comment}

\usepackage{listings}

\newcommand{\Sym}{\textnormal{Sym}}

\newcommand{\eps}{\varepsilon}

\newcommand{\inParens}[1]{{\left({#1}\right)}}
\newcommand{\inBrackets}[1]{{\left[{#1}\right]}}
\newcommand{\inBraces}[1]{{\left\{{#1}\right\}}}
\newcommand{\inAbs}[1]{{\left\lvert{#1}\right\rvert}}

\newcommand{\tucker}{\bullet}
\newcommand{\flatten}[2]{{#1}\times{#2}}

\newcommand{\todo}[2]{{\color{blue}\textnormal{\textbf{TODO}({#1}): {#2}}}}

\newcommand{\pa}{\textnormal{pa}}

\newcommand{\nd}{\textnormal{nd}}
\renewcommand{\top}{\textnormal{top}}

\newcommand{\RtoTheV}{{\bb R^{\inAbs{V}}}}

\begin{document}
\maketitle
\begin{abstract}In this paper we study the space of second- and third-order moment tensors of random vectors which satisfy a Linear Non-Gaussian Acyclic Model (LiNGAM). In such a causal model each entry $X_i$ of the random vector $X$ corresponds to a vertex $i$ of a directed acyclic graph $G$ and can be expressed as a linear combination of its direct causes $\{X_j: j\to i\}$ and random noise. 
For any directed acyclic graph $G$, we show that a random vector $X$ arises from a LiNGAM with graph $G$ if and only if certain easy-to-construct matrices, whose entries are second- and third-order moments of $X$,
drop rank. 
This determinantal characterization extends previous results proven for polytrees  and generalizes the well-known local Markov property for Gaussian models.
\end{abstract}

\section{Introduction}

Structural equation models (SEMs)~\cite{maathuis2018handbook} capture cause-effect relationships among a set of random variables $\{X_i, i\in V\}$ by hypothesizing that each variable
is a noisy function of its direct causes. Given a directed acyclic graph (DAG) $G = (V, E)$ with one random variable $X_i$ associated to each vertex $i\in V$, the directed edges $j\to i\in E$ correspond to direct causes. A {\em linear structural equation model} with graph $G$ then hypothesizes that
\begin{align}\label{eq:structural_equations}
X_i = \sum_{j\to i\in E}\lambda_{ji}X_j + \varepsilon_i,\quad i\in V,
\end{align}
where the random variables $\varepsilon_i$ are mutually independent and represent random noise.

Classically, the noise terms $\varepsilon_i$ are assumed to be Gaussian, in which case one aims to learn the graph $G$ solely from the covariance matrix of $X$. In this scenario, one can only learn $G$ up to a Markov equivalence class. On the other hand, linear {\em non-Gaussian} acyclic models (LiNGAM)~\cite{shimizu2006linear} assume that the error terms $\varepsilon_i$ are non-Gaussian. They have sparked a wide range of interest since they allow one to learn the true directed acyclic graph $G$ rather than its Markov equivalence class from observational data~\cite{shimizu2006linear}.


One way of learning the graph $G$ both in the Gaussian and non-Gaussian settings is via the method of moments.
In such methods, one obtains insights about the algebraic structure of the moments of the random vector $X$ for each graph~\cite{drton2017structure, sullivant2018algebraic, wang2019highdimensional, wang:2023, schkoda2023goodnessoffit}, and then devises an algorithm that utilizes these insights and learns the graph.  
The algebraic relations that hold among the entries of the covariance matrix of $X$ have been a major topic of interest in the algebraic statistics community ~\cite{drton2020nested,sullivant2008algebraicGeometry,vanOmmen2017algebraic, draisma2013positivity, sullivant2010trek}. 

In this work we focus on the set of second- and third-order moments of a random vector $X$ which arises from a LiNGAM for a given DAG $G$. These moments are enough to identify the graph $G$
and to learn it efficiently even in the high-dimensional setting
~\cite{wang2019highdimensional}. Specific low-degree determinantal relationships among the second and third order moments for a linear non-Gaussian causal model have given rise to efficient algorithms for learning the graph $G$ from data~\cite{schkoda2024causal, wang2019highdimensional, wang:2023, drton2025cyclic}, but a complete algebraic characterization of the model of second- and third-order moments is only known when the graph $G$ is a polytree~\cite{amendola2023third}. 

We here complete this characterization for any DAG $G$. 
We show that rank constraints on certain matrices whose entries are second- or third-order moments of the random vector $X$ uniquely specify the DAG $G$ (see Theorem~\ref{thm:main}). Our constraints contain as a subset the well-known constraints arising from the local Markov property satisfied by the covariance matrix only~\cite{sullivant2018algebraic}, and also extend recent work on LiNGAM where $G$ is assumed to be a polytree~\cite{amendola2023third}.

We first illustrate our result in an example.

\begin{example}\label{ex:1}

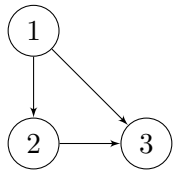
\begin{figure}[h]
    \centering
    \begin{tikzpicture}
        \tikzset{vertex/.style = {shape=circle,draw,minimum size=1.5em}}
        \tikzset{edge/.style = {->,> = latex'}}
        \node[vertex] (1) at  (  0,1.5) {$1$};
        \node[vertex] (2) at  (  0,  0) {$2$};
        \node[vertex] (3) at  (1.5,  0) {$3$};
        \draw[edge] (1) to (2);
        \draw[edge] (1) to (3);
        \draw[edge] (2) to (3);
    \end{tikzpicture}
    \caption{
        The complete DAG on $3$ vertices.
    }
    \label{fig 3 node acyclic}
\end{figure}

Consider the graph $G$ with vertices $V = \{1,2,3\}$ and edges $\{1\to2, 1\to 3, 2\to 3\}$, the simplest DAG which is not a polytree (see Figure~\ref{fig 3 node acyclic}). 
Our Theorem~\ref{thm:main} implies that if $X$ lies in the LiNGAM for some DAG, then it lies in the LiNGAM for this particular DAG $G$ if and only if its
set of second-order moments $s_{ij} = \mathbb E[X_iX_j]$ and third-order moments $t_{ijk} = \mathbb E[X_iX_jX_k]$ are such that the matrices
$$M_2 = \begin{bmatrix}s_{11} & t_{111} & t_{112} & t_{113}\\
s_{12} & t_{112} & t_{122} & t_{123}
\end{bmatrix}, \,\, M_3 = \begin{bmatrix}
s_{11}&s_{12} & t_{111}&t_{112}&t_{113}&t_{122}&t_{123}\\
s_{12}&s_{22}&t_{112}&t_{122}&t_{123}&t_{222}&t_{223}\\
s_{13}&s_{23}& t_{113}&t_{123}&t_{133}&t_{223}&t_{233}
\end{bmatrix}$$
drop rank, i.e., they have ranks $1$ and $2$, respectively. 
\end{example}

 The rest of this paper is organized as follows. We begin in Section~\ref{sec:preliminaries} by a description of linear non-Gaussian acyclic models as well as the parametrization they imply for the second- and third-order moments of the random vector $X$. We include a review of relevant prior work on linear Gaussian or 
 non-Gaussian models. In Section~\ref{sec:results}, we present our main result, Theorem~\ref{thm:main}, which exhibits the constraints that characterize the set of second- and third-order moments corresponding to a given DAG $G$. We note that our result includes the constraints arising from the local Markov property in the Gaussian case~\cite{sullivant2018algebraic}. In Section~\ref{sec:additional_equations}, we show how to derive additional polynomial equations from the ones implied by Theorem~\ref{thm:main}. In particular, we show how our result generalizes the characterization of the defining equations for polytrees in the non-Gaussian case~\cite{amendola2023third}.
 While previous work~\cite{wang2019highdimensional} uses algebraic constraints to recover the sources of a DAG, in Section~\ref{sec:application} we show how to use our result for the recovery of sink nodes.
 We conclude in Section~\ref{sec:discussion} with a discussion and further questions of interest.
 The proofs of our results are located in Section~\ref{sec:proofs}.
 
\section{Background}\label{sec:preliminaries}
In this section we first introduce the mathematical formulation of our problem and prior work related to it.
\subsection{Preliminaries}
 
Let $G = (V, E)$ be a DAG, and let $(X_i, i \in V)$, be a collection of random variables indexed by the vertices in $V$. A vertex $j \in V$ is a {\em parent} of a vertex $i$ if there
is an edge pointing from $j$ to $i$, i.e., if $(j, i) \in E$ which we will also write as $j \to i \in E$. A vertex $j\in V$ is a {\em non-descendant} of a vertex $i$ if there are no directed paths from $i$ to $j$. We denote
the set of all parents of vertex $i$ by $\text{pa}(i)$ and the set of all non-descendants of $i$ by $\text{nd}(i)$.  As described in~\eqref{eq:structural_equations}, the graph $G$ gives rise to the linear structural equation model
consisting of the joint distributions of all random vectors $X = (X_i, i \in V )$ such that
\begin{align*}
X_i = \sum_{j\in\text{pa}(i)}\lambda_{ji}X_j + \varepsilon_i, \,\, i\in V,
\end{align*}
where the $\varepsilon_i$ are mutually independent random variables representing stochastic errors. The errors
are assumed to have expectation $\mathbb E[\varepsilon_i] = 0$, finite variance $\omega^{(2)}_i = \mathbb E[\varepsilon_i^2] > 0$, and finite third
moment $\omega^{(3)}_i = \mathbb E[\varepsilon_i^3]$. No other assumption about their distribution is made and, in particular,
the errors need not be Gaussian (in which case we would have $\mathbb E[\varepsilon_i^3 ] = 0$ by symmetry of the
Gaussian distribution). The coefficients $\lambda_{ji}$ in~\eqref{eq:structural_equations} are unknown real-valued parameters, and we
fill them into a matrix $\Lambda = (\lambda_{ji})\in\mathbb R^{|V|\times |V|}$ by adding a zero entry when $(j, i) \not\in E$. We denote the
set of all such sparse matrices as $\mathbb R^E$ . We note that for simplicity, and without loss of generality,
the equations in~\eqref{eq:structural_equations} do not include a constant term, so we have $\mathbb E[X_i] = 0$ for all $i \in V$.

The structural equations~\eqref{eq:structural_equations} can be rewritten in matrix-vector format as $X = \Lambda^TX+\varepsilon$, and, thus,
$$X = (I-\Lambda)^{-T}\varepsilon,$$
where we note that the matrix $I-\Lambda$ is always invertible when the graph $G$ is acyclic~\cite{sullivant2018algebraic}.

Let $\Omega^{(2)}= (\mathbb E[\varepsilon_i\varepsilon_j])$ and $\Omega^{(3)} = (\mathbb E[\varepsilon_i\varepsilon_j \varepsilon_k])$ be the covariance matrix and the tensor of third-order
moments of $\varepsilon$, respectively. Both $\Omega^{(2)}$ and $\Omega^{(3)}$ are diagonal, with the diagonal entries being $\Omega^{(2)}_{ii} = \omega^{(2)}_i = \mathbb E[\varepsilon_i^2] >0$ and $\Omega^{(3)}_{iii} = \omega^{(3)}_i = \mathbb E[\varepsilon_i^3]$.
\begin{lemma}\label{lem:parametrization}
The covariance matrix and the third-order moment tensor of the solution $X$ of~\eqref{eq:structural_equations}
are equal to
\begin{align*}
S &= (s_{ij}) = (I - \Lambda)^{-T} \Omega^{(2)}(I - \Lambda)^{-1},\\
T &= (t_{ijk}) = \Omega^{(3)} \bullet (I - \Lambda)^{-1} \bullet (I - \Lambda)^{-1} \bullet (I - \Lambda)^{-1},
\end{align*}
respectively. Here $\bullet$ denotes the Tucker product~\cite{kolda2009tensor}.
\end{lemma}
This fact follows from standard results on how moments of random vectors change after
linear transformation. For a complete proof, see, e.g.,~\cite[Proposition 1.2]{amendola2023third}.

As we are assuming positive error variances, $\mathbb E[\varepsilon^2_i ] > 0$, the matrix $\Omega^{(2)}$ is
positive definite and the same is true for the covariance matrix $S$ of $X$. Since $\Omega^{(3)}$ is diagonal, the
third-order moment tensor $T$ of $X$ is a symmetric tensor of symmetric tensor rank at most $\inAbs{V}$; this need
not be the case for a general $|V| \times |V| \times |V|$ tensor~\cite{comon2008symmetric}. In the sequel, we write $PD(\RtoTheV)$ for the positive
definite cone in $\mathbb R^{|V|\times |V|}$ and $\text{Sym}^3(\RtoTheV)$ for the space of symmetric tensors in $\mathbb R^{|V| \times |V| \times |V|}$.

\begin{definition}
Let $G = (V, E)$ be a DAG. The {\em second- and third-order moment model of $G$} is the set $\mathcal M^{\leq 3}(G)$
that comprises of all pairs of covariance matrices $S$ and third-order moment tensors $T$ that are realizable under
the linear structural equation model given by $G$. That is,
$$\mathcal M^{\leq 3}(G) = \{(\underbrace{(I -\Lambda)^{-T}\Omega^{(2)}(I - \Lambda)^{-1}}_{S}, \underbrace{\Omega^{(3)} \bullet (I - \Lambda)^{-1} \bullet (I - \Lambda)^{-1} \bullet (I - \Lambda)^{-1}}_{T}) :$$
$$\Omega^{(2)}\in PD(\RtoTheV) \,\, \text{diagonal}, \Omega^{(3)} \in \text{Sym}^3(\RtoTheV)\,\, \text{diagonal}, \Lambda \in \mathbb R^E \}$$$$ \subseteq PD(\RtoTheV) \times \text{Sym}^3(\RtoTheV).$$
Furthermore, the second- and third-order moment ideal of $G$ is the ideal $\mathcal I^{\leq 3}(G)$ of polynomials in the entries
$S = (s_{ij} )$ and $T = (t_{ijk} )$ that vanish when $(S, T )\in\mathcal M^{\leq 3}(G)$.
\end{definition}

The problem which we solve here is as follows.

\begin{problem}
Assume that $S\in PD(\RtoTheV)$ is a positive definite matrix and $T\in \Sym^3(\RtoTheV)$ is a symmetric tensor. 
Given a DAG $G$, find polynomial constraints in the entries of $S$ and $T$ which 
are satisfied if and only if $(S, T)$ lies in $\mathcal M^{\leq 3}(G)$.
\end{problem}

\subsection{Prior work}
Here we summarize relevant prior work on the algebraic description of the model of interest $\mathcal M^{\leq 3}(G)$.
\subsubsection*{Linear Gaussian models}
This problem has classically been studied in the algebraic statistics literature in the case of Gaussian graphical models~\cite{sullivant2018algebraic}. Here, they consider a linear structural equation model with Gaussian error terms $\varepsilon_i$. As a result, the third-order moments all vanish, and the problem is to describe the model $\mathcal M^{2}(G)$ consisting of all covariance matrices $S$ which factorize as $S = (I -\Lambda)^{-T}\Omega^{(2)}(I - \Lambda)^{-1}$ with $\Omega^{(2)}\in PD(\RtoTheV)$ and $\Lambda\in\mathbb R^E$.

Conditional independence implies constraints on the entries of $S$ for a given graph $G$ as follows. If sets of vertices $A$ and $B$ are $d$-separated given set $C$ in the graph $G$ (see, e.g.~\cite{maathuis2018handbook} for the definition of $d$-separation), then the Global Markov Property implies that $X_A$ is conditionally independent of $X_B$ given $X_C$. For a Gaussian distribution this is equivalent to the submatrix $S_{A\cup C, B\cup C}$ of the covariance matrix $S$ with rows indexed by $A\cup C$ and columns indexed by $B\cup C$ having rank at most the size of $C$, $\inAbs{C}$, as shown in~\cite{sullivant2008algebraicGeometry}. Furthermore, the model $\mathcal M^2(G)$ is cut out by precisely these rank constraints arising from all the different $d$-separation statements which hold for $G$ inside the positive definite cone $PD(\RtoTheV)$~\cite{sullivant2008algebraicGeometry}.

The seminal paper~\cite{sullivant2010trek} then asks the question of what other rank constraints hold for submatrices of the covariance matrix $S$ in a Gaussian graphical model, and the answer comes via trek-separation.

\begin{definition}[\cite{sullivant2010trek, robeva2021multi-trek}]
Given $k\geq 2$ vertices $v_1,\ldots v_k$, a {\em $k$-trek} $\mathcal T$ between them is an ordered tuple of directed paths $(P_1,\ldots, P_k)$ which have a common source node $t$
, called the top of $\cal T$ ($\top(\cal T)$), and the path $P_i$ has sink $v_i$ for each $i = 1,\dots, k$. 
A 2-trek is usually known as a {\em trek}.

Let $A, B\subseteq V$ be two subsets of vertices. The pair of sets $(L, R)$ {\em trek-separates} $A$ and $B$
 if for every trek $\mathcal T = (P_1,P_2)$ between a vertex $a\in A$ and a vertex $b\in B$ either $P_1$ contains a vertex from $L$ or $P_2$ contains a vertex from $R$. \end{definition}
The main result in~\cite{sullivant2010trek} states that the submatrix $S_{A,B}$ has rank at most $r$ if and only if there exist sets of vertices $L, R\subseteq V$ such that $(L, R)$ trek-separates $(A,B)$ and $\inAbs{L} + \inAbs{R}\leq r$. Therefore, all rank constraints on $S$ correspond to trek-separation in the graph.

When the graph $G$ has hidden variables, however, rank constraints are not enough to cut out the model $\mathcal M^2(G)$~\cite{sullivant2008algebraicGeometry}. Such constraints, also knwon as Verma constraints, can sometimes be expressed in the form of nested determinants~\cite{drton2020nested}, but a general way of obtaining all of them is not known.

\subsubsection*{Linear Non-Gaussian Models}
Such a thorough study has not been done for the LiNGAM models $\mathcal M^{\leq 3}(G)$ which we consider here. Instead of only looking at the covariance matrix $S$, we now also have access to the third-order moment tensor $T$. The constraints which cut out the model $\mathcal M^{\leq 3}(G)$ have been discovered in the case when $G$ is a polytree~\cite{amendola2023third}. For arbitrary DAGs $G$,  a generalized version of trek-separation has been found~\cite{robeva2021multi-trek}, but a complete characterization of the constraints that cut out $\mathcal M^{\leq 3}(G)$ was not known prior to the present work.

On the algorithmic side, methods based on algebraic constraints which hold among the second- and higher-moments of the random vector $X$ have found success. The work~\cite{wang2019highdimensional} develops a high-dimensional algorithm for learning the DAG $G$ based on testing the rank of certain $2\times 2$ matrices (see Corollary~\ref{cor:sources}) in order to successively find source nodes and remove them from the graph. More recently,~\cite{schkoda2024causal} uses rank constraints on matrices which consist of second-, third-, and higher-order moments in order to learn a DAG $G$ with hidden variables. Furthermore, ~\cite{schkoda2023goodnessoffit} uses algebraic constraints for goodness-of-fit tests which determine whether the data arises from a linear non-Gaussian model at all.









\subsection{Notation}
We will be interested in matrices whose entries consist of blocks of $S$, and blocks of $T$. 
Thus, we require notation for blocking out sections of matrices and tensors.
For two subsets $A,B \subseteq V$, we can define $S_{A,B}$
to be the matrix with entries $S_{a,b}$ with $a \in A$, and $b \in B$. Here
$a$ is the row label, and $b$ is the column label.

We will need to do a similar operation with $T$.
For three subsets $A,B,C \subseteq V$, we define $T_{A,\flatten{B}{C}}$
to be a matrix with entries $T_{a,b,c}$ where $a \in A$, $b \in B$, and $c \in C$ 
which is flattened such that $a$ is the row label and $(b,c)$ is the column label. 
As an example, let $A = B = C = \{1, 2\}$.
Then $T_{A,\flatten{B}{C}}$ is
\begin{equation*}
    T_{\{1, 2\},\flatten{\{1, 2\}}{\{1, 2\}}}
    =
    \begin{blockarray}{ccccc}
        &(1,1) &(1,2) &(2,1) &(2,2)
        \\
        \begin{block}{c[cccc]}
            1 &T_{1,1,1} &T_{1,1,2} &T_{1,2,1} &T_{1,2,2} 
            \\
            2 &T_{2,1,1} &T_{2,1,2} &T_{2,2,1} &T_{2,2,2} 
            \\
        \end{block}
    \end{blockarray}.
\end{equation*}

\section{Algebraic characterization of $\mathcal M^{\leq 3}(G)$}\label{sec:results}

We are now ready for our main result. The following theorem gives explicit constraints that cut out the set $\mathcal M^{\leq 3}(G)$ of pairs $(S, T)$ that arise from the LiNGAM with DAG $G$. 
\begin{theorem}
    Let $G = (V,E)$ be a DAG, $S \in PD(\RtoTheV)$, and $T \in \Sym^3(\RtoTheV)$.
    Then, $(S, T)$ lies in the model $\cal M^{\leq 3}(G)$ if and only if
    for every vertex $v \in V$, the following matrix has rank equal to $|\pa(v)|$,
    \begin{equation}
        M_v := \begin{bmatrix}
                S_{\pa(v),\nd(v)} &T_{\pa(v),\flatten{\nd(v)}{V}} 
                \\
                S_{v,\nd(v)} &T_{v,\flatten{\nd(v)}{V}}
            \end{bmatrix}.
        \label{equ definition of M_v}
    \end{equation}
    \label{thm:main}
\end{theorem}
\begin{remark}
    The rank constraints given in~\eqref{equ definition of M_v} 
    are equivalent to the vanishing of all $(\inAbs{\pa(v)}+1)$-minors of the matrix $M_v$. 
    Since $S$ is assumed to be positive definite, 
    the sub-matrix $S_{\pa(v),\pa(v)}$ is invertible, 
    and therefore the rank constraints are also equivalent to the last row of $M_v$ being a linear combination of its other rows.
    \label{rmk:postTheorem}
\end{remark}

\begin{example}(Example~\ref{ex:1} Continued)
Consider the graph $G$ with vertices $V = \{1,2,3\}$ and edges $\{1\to2, 1\to 3, 2\to 3\}$ (Figure~\ref{fig 3 node acyclic})). 
Using Macaulay2, we can define the ideal $\mathcal I^{\leq 3}(G)$ which equals the kernel of the parametrization of the moments $S$ and $T$ in terms of the parameters $\Lambda$, $\Omega^{(2)}$, and $\Omega^{(3)}$  from Lemma~\ref{lem:parametrization}. We can then confirm that this ideal equals the ideal generated by the $2-$minors of the matrix $M_2$ and the $3-$minors of the matrix $M_3$ that arise from our Theorem~\ref{thm:main}
$$M_2 = \begin{bmatrix}s_{11} & t_{111} & t_{112} & t_{113}\\
s_{12} & t_{112} & t_{122} & t_{123}
\end{bmatrix}, \,\, M_3 = \begin{bmatrix}
s_{11}&s_{12} & t_{111}&t_{112}&t_{113}&t_{122}&t_{123}\\
s_{12}&s_{22}&t_{112}&t_{122}&t_{123}&t_{222}&t_{223}\\
s_{13}&s_{23}& t_{113}&t_{123}&t_{133}&t_{223}&t_{233}
\end{bmatrix},$$
saturated by the principal minors of the matrix $S$. 
Indeed, according to Theorem~\ref{thm:main}, the model $\mathcal M^{\leq 3}(G)$ is cut out by the minors of these matrices inside $PD(\RtoTheV)\times\text{Sym}^3(\RtoTheV)$, and so the principal minors of $S$ are always nonzero.
\end{example}

\begin{remark}
    \newcommand{\sep}{\textnormal{sep}}
    Fix $v \in V$.
    The block of $M_v$ containing only entries of $S$
    is as follows,
    \begin{equation*}
        \begin{blockarray}{ccc}
            &\pa(v) &\nd(v)\setminus\pa(v) \\
            \begin{block}{c[cc]}
                \pa(v) &S_{\pa(v),\pa(v)} &S_{\pa(v),\nd(v)\setminus\pa(v)}
                \\
                v       &S_{v,\pa(v)}     &S_{v,\nd(v)\setminus\pa(v)}
                \\
            \end{block}
        \end{blockarray}.
    \end{equation*}
    The fact that it drops rank represents exactly the Local Markov Property which implies that 
    $X_v$ is conditionally independent from $X_{\nd(v)\setminus\pa(v)}$ given $X_{\pa(v)}$~\cite{maathuis2018handbook,sullivant2018algebraic}.
\end{remark}

The following Lemma is the necessity condition of Theorem~\ref{thm:main}.
\begin{lemma}
    Let $G = (V,E)$ be a DAG, if $(S, T)$ lies in the model $\cal M^{\leq 3}(G)$, then
    for any vertex $v \in V$, the following matrix has rank $\inAbs{\pa(v)}$,
    \begin{equation}
        M_v := \begin{bmatrix}
                S_{\pa(v),\nd(v)} &T_{\pa(v),\flatten{\nd(v)}{V}} 
                \\
                S_{v,\nd(v)} &T_{v,\flatten{\nd(v)}{V}}
            \end{bmatrix}.
    \end{equation}
    Since these matrices drop rank, all $(\inAbs{\pa(v)}+1) \times (\inAbs{\pa(v)}+1)$ minors must vanish.
    \label{lemma matrix with low rank}
\end{lemma}
\begin{proof}
    Recall that $X_v = \sum_{u\in\pa(v)}\lambda_{uv} X_u + \eps_v$,
    If $w$ is any non-descendant of $v$, $X_w$ is independent of $\eps_v$, and so $\bb E[\eps_v X_w ] = 0$.
    We see that
    \begin{equation*}
        s_{vw} = \mathbb E[X_vX_w] =  \sum_{u\in\pa(v)}\lambda_{uv} \mathbb E[X_uX_w] + \mathbb E[\eps_vX_w]
        = \sum_{u\in\pa(v)} \lambda_{uv} s_{uw}.
    \end{equation*}
    A similar equation can be derived for $t_{vwz}$ for any $z\in V$:
    $$t_{vwz} = \sum_{u\in\pa(v)}\lambda_{uv}t_{uwz}.$$
Then, the vector $(-1, \lambda_{\pa(v), v})^T$ is in the left null space of the matrix $M_v$, which means that $M_v$ has rank at most $|\pa(v)|$. 

In a DAG, $\pa(v) \subseteq \nd(v)$,
    and so $S_{\pa(v),\pa(v)}$ is a submatrix of $M_v$, which does not drop rank since $S$ is positive definite.
    Therefore $M_v$ 
    is a $\inParens{\inAbs{\pa(v)} + 1}$ by $\inAbs{\nd(v)}\inParens{\inAbs{V}+1}$ matrix 
    of rank exactly equal to $\inAbs{\pa(v)}$.
\end{proof}

The next Lemma is the more difficult step in the proof of Theorem~\ref{thm:main}.
We show that if the matrices in equation~\eqref{equ definition of M_v} drop rank,
then we can construct the appropriate $\Lambda, \Omega^{(2)}, \Omega^{(3)}$ which parametrize $S$ and $T$.
This guarantees that $(S,T)$ will lie in the model $\cal M^{\leq 3}(G)$.

\begin{lemma}\label{lem:main}
    Let $(G, E)$ be a directed acyclic graph.
    Let $(S, T) \in PD(\RtoTheV) \times \Sym_3(\RtoTheV)$.
    For each vertex $v \in V$, let $M_v$ be the matrix defined 
    in equation~\eqref{equ definition of M_v}.
    Suppose that each $M_v$ has rank $\inAbs{\pa(v)}$. 
    Then, there exists $\Lambda \in \bb R^E$, such that
    \begin{align}
        \Omega^\inParens{2} 
        &:= (I-\Lambda)^TS (I-\Lambda),
        \label{equ diagonalization of S}
        \\
        \Omega^\inParens{3} 
        &:= T \tucker (I-\Lambda) \tucker (I-\Lambda) \tucker (I-\Lambda)
        \label{equ diagonalization of T}
    \end{align}
    are diagonal.
    \label{lem matrices recover the graph}
\end{lemma}
\begin{proof}[Proof sketch]
    Fix $v \in V$.
    By Remark~\ref{rmk:postTheorem}, 
    the bottom row of the matrix $M_v$ is a linear combination of the other rows. 
    We define $\lambda_{iv}$ to be the coefficient of the $i$-th row in this linear combination,
    and all other $\lambda_{kv} = 0$. 
    In this way, we define the matrix $\Lambda\in\bb R^E$.
    A calculation then shows that $\Omega^\inParens{2}$ and $\Omega^\inParens{3}$ are both diagonal.
\end{proof}

Theorem~\ref{thm:main} is established by combining Lemma~\ref{lemma matrix with low rank} and Lemma~\ref{lem matrices recover the graph}. The full proof are located in Section~\ref{sec:proofs}.

\subsection*{Computing the vanishing ideal of $\mathcal M^{\leq 3}(G)$}
While in Theorem~\ref{thm:main} we found algebraic constraints that cut out our model $\mathcal M^{\leq 3}(G)$, we have not discussed its vanishing ideal. We conjecture that the vanishing ideal of the model $\mathcal M^{\leq 3}(G)$ equals the ideal generated by the minors of the matrices $M_v$ for $v\in V$ saturated by the product of all principal minors of the matrix $S$.

The paper~\cite{boege2024realbirationalimplicitizationstatistical} describes a principled way of finding this ideal exactly assuming the parameters $\Lambda, \Omega^{(2)}, \Omega^{(3)}$ are identifiable from $(S, T)$. This is indeed the case for us, and Theorem 3.11 of~\cite{boege2024realbirationalimplicitizationstatistical} can be applied. (We would follow their Section~4.2, augmented by the variables $t_{ijk}$. To express our parametrization as a birational map, we include all entries of $\Omega^{(3)}$ in the domain.) However, the equations that one obtains using Theorem 3.11 from~\cite{boege2024realbirationalimplicitizationstatistical} appear more complicated than the principal minors of our matrices $M_v$. Therefore, we leave it as future work to compute the vanishing ideal of the model $\mathcal M^{\leq3}(G)$.

\section{Additional equations}\label{sec:additional_equations}
While we have derived enough equations in Theorem~\ref{thm:main} to cut out our model $\mathcal M^{\leq 3}(G)$, in this section we show how to find additional equations which also vanish on our model, but are not  minors of the matrices $M_v$ defined in Theorem~\ref{thm:main}.
\begin{definition}
    For $A, B \subseteq V$, define
    \begin{equation*}
        R_{A,B} := \begin{bmatrix}
            S_{A,B} &T_{A,B\times V}
        \end{bmatrix}.
    \end{equation*}
\end{definition}
In the statement of Theorem~\ref{thm:main}, we have matrices of the form
\begin{equation*}
    M_v =
    \begin{bmatrix}
        R_{\pa(v),\nd(v)}
        \\
        R_{v,\nd(v)}
    \end{bmatrix}.
\end{equation*}
The paper~\cite{amendola2023third} shows that the model $\mathcal M^{\leq 3}(G)$ corresponding to a polytree $G$ is cut out by only $2\times 2$ determinants inside $PD(\RtoTheV)\times \Sym^3(\RtoTheV)$. However, the rank conditions from Theorem~\ref{thm:main} may involve larger determinants even in the case of polytrees.
Thus, we here ask if we can replace $\pa(v)$, and $\nd(v)$ to obtain different matrices whose minors also vanish on the model.

\begin{proposition}
    For any $A \subseteq V$, define
    \begin{equation*}
        \pa(A) := \inParens{\bigcup_{a\in A}\pa(a)}\setminus A = \inBraces{\textit{all vertices outside $A$ with an edge pointing into $A$}}.
    \end{equation*}
    For any $A \subseteq V$, define
    \begin{equation*}
        \nd(A) := \bigcap_{a\in A}\nd(a) =  \inBraces{\textit{all common non-descendants of the vertices in $A$}}.
    \end{equation*}
    Let $v \in G$, and $A \subseteq V \setminus  v$.
    Then the last row of the following matrix is a linear combination of the other rows.
    \begin{equation*}
        \begin{bmatrix}
                R_{\pa(\{v\}\cup A),\nd(\{v\}\cup A)}
                \\
                R_{v,\nd(\{v\}\cup A)}
        \end{bmatrix}.
    \end{equation*}

    \label{prop other matrices}
\end{proposition}
\begin{remark}
    Given a collection of vertices $A$, 
    we can think of the above equations as corresponding to removing the set $A$ from the structural equation model.
    As an example, let $v \in V$, and let $a\in V$ be a parent of $v$. 
    Then
    \begin{equation*}
        X_v = \lambda_{av} X_a + \dots + \eps_v.
    \end{equation*}
    Removing the vertex $a$ from the structural equation model 
    would require replacing $X_a$ with its expression involving its parents, 
    and augmenting $\eps_v$ to $\eps_v + \lambda_{av} \eps_a$.
    In this way, we have expanded $\pa(v)$ to $\pa(\inBraces{v,a})$. Since $\eps_v + \lambda_{av} \eps_a$ is only independent of those $X_w$ for which $w$ is a non-descendant of both $v$ and $a$, we are forced to shrink $\nd(v)$ to $\nd(\inBraces{v,a})$. For more details, see the proof sketch of Lemma~\ref{lemma matrix with low rank} located in Section~\ref{sec:results}.
\end{remark}

Proposition~\ref{prop other matrices} can be applied to generate more equations.
For example, if $w \in V$ is a source vertex, we might ask if there is a collection of equations which determine this fact.

\begin{corollary}[Sources]\label{cor:sources}
    Let $w \in V$ be a source of the graph $G$, $v \in V\setminus w$.
    Then the following matrix has rank $1$,
    \begin{equation*}
        \begin{bmatrix}
              s_{w,w} &t_{w,w,w}
            \\s_{v,w} &t_{v,w,w}
        \end{bmatrix}.
    \end{equation*}
\end{corollary}
\begin{proof}
    We set $A = V\setminus \{v,w\}$.
    Then, $\{w\} \supseteq \pa(\{v\}\cup A)$,
    and $w \in \nd(\{v\}\cup A)$.
    Applying Proposition~\ref{prop other matrices} gives the result.
\end{proof}

 This is the precise equation used in the algorithm proposed in~\cite{wang2019highdimensional}, which efficiently recovers the DAG $G$ by recursively identifying source nodes.

\subsection*{Polytrees}

A polytree $G$ is a directed graph whose undirected skeleton is a tree. Such a graph is naturally acyclic.
Polytrees also have the property that for any two $v,w \in V$, there is at most one simple trek (a trek whose left and right side do not share any edges) between $v$ and $w$. We denote the top of this trek by $\top(v,w)$. The ideal which cuts out the model $\mathcal M^{\leq 3}(G)$ inside $PD(\mathbb R^{|V|})\times \Sym^3(\mathbb R^{|V|})$ of a polytree $G$ is given in~\cite{amendola2023third}.

\begin{proposition}[{\citealp[Lemma~3.4~(b)]{amendola2023third}}]
    Let $G = (V, E)$ be a polytree, and $v, w \in V$. 
    Then if there is an edge between $v$ and $w$, the $2$-minors of the following trek-matrix vanish,
    \begin{equation*}
        \begin{bmatrix}
               s_{vk_1} &\dots &s_{vk_r} &t_{v\ell_1 m_1} &\dots &t_{v\ell_q m_q}
            \\ s_{wk_1} &\dots &s_{wk_r} &t_{w\ell_1 m_1} &\dots &t_{w\ell_q m_q}
        \end{bmatrix},
    \end{equation*}
    where
    \begin{itemize}
        \item $k_1, \dots, k_r$ are vertices such that $\top(v, k_a) = \top(w, k_a)$ for $a = 1, \dots , r$, and
        \item $(\ell_1, m_1), \dots, (\ell_q , m_q)$ are such that $\top(v, \ell_b, m_b) = \top(w, \ell_b, m_b)$ for $b = 1, \dots, q$.
    \end{itemize}
    \label{prop:polyTreeEquations}
\end{proposition}


We show how to recover the result of Proposition~\ref{prop:polyTreeEquations} using Proposition~\ref{prop other matrices}.
\begin{proof}
    Let $v,w \in V$ such that there is an edge between them.
    Without loss of generality, assume that $w$ is a parent of $v$.
    Let $A \subset V$ be all $a \in V\setminus \{v\}$ such that there is a directed path from $a$ to $v$ which does not pass through $w$.
    We claim that $\pa(\{v\}\cup A) = \{w\}$.

    Since $w$ is a parent of $v$ and $w\not\in A$, then $w \in \pa(\{v\}\cup A)$.
    Now for any $u \in\pa(\{v\}\cup A)$, we know that $u \notin \{v\}\cup A$, and there is an edge from $u$ to either $a=v$ or any element $a$ of $A$.
    By definition of $A$, there is a directed path $\gamma$ from $a$ to $v$ which does not go though $w$.
    Thus, $u \to \gamma$ is a directed path from $u$ to $v$ which can only pass through $w$ if $u = w$.
    Since $u \notin \{v\}\cup A$, then $u = w$.

    We also claim that $\nd(\{v\}\cup A)$ contains all $u \in V$ such that $\top(w,u) = \top(v,u)$.
    Let $u$ be such a vertex.
    If $u$ is a descendant of $v$, then $\top(v,u) = v$ but $\top(w,u) = w$.
    Therefore, $u$ is not a descendant of $v$.
    If there exists $a\in A$ such that $u$ is a descendant of $a$, then $\top(v,u) \in A$.
    Since the skeleton of $G$ is a tree, there is only one (undirected) path between $w$ and $u$, which must go through $v$.
    Thus, there is no trek between $w$ and $u$, and so there is no $\top(w,u)$.
    Hence $u \in \nd(\{v\}\cup A)$.
    
    Applying Proposition~\ref{prop other matrices} shows that the following matrix has rank at most 1,
    \begin{equation}
        \begin{bmatrix}
              S_{w,\nd(\{v\}\cup A)} &T_{w,\nd(\{v\}\cup A)\times V}
            \\S_{v,\nd(\{v\}\cup A)} &T_{v,\nd(\{v\}\cup A)\times V}
        \end{bmatrix}.
        \label{equ:polytreesMatrix}
    \end{equation}
    By the arguments above, the matrix in Proposition~\ref{prop:polyTreeEquations} is a submatrix of the matrix in equation~\eqref{equ:polytreesMatrix} with the same number of rows and possibly fewer columns. Therefore, it has rank at most 1 as well.
\end{proof}

\section{Proofs}\label{sec:proofs}
In this section we prove our main results Theorem~\ref{thm:main} and Proposition~\ref{prop other matrices}.
\subsection{Proof of Theorem~\ref{thm:main}}

We can now finish the proof of Theorem~\ref{thm:main}. 
Lemma~\ref{lem matrices recover the graph} shows that if the matrices $M_v$ all drop rank, 
then we can parametrize $(S,T)$ in the way which guarantees that $(S,T) \in \cal M^{\leq 3}(G)$. 
\begin{proof}[Proof of Lemma~\ref{lem matrices recover the graph}]
    In a directed acyclic graph, $\pa(v) \subseteq \nd(v)$.
    Thus, $S_{\pa(v),\pa(v)}$ is a submatrix of $M_v$. Since $S_{\pa(v), \pa(v)}$ is invertible, then the bottom row of $M_v$ is a linear combination of the top $|\pa(v)|$ rows. Let $\lambda_{iv}$ be the coefficients of this linear combination. Set the remaining $\lambda_{kv}$ to $0$ and define $\Lambda$ to be the matrix containing these coefficients. Then $\Lambda \in \bb R^E$.
    
    

    Now, let $\Omega^{(2)} = (I-\Lambda)^TS(I-\Lambda)$ as in equation~\eqref{equ diagonalization of S}. 
    We will show that $\Omega^{(2)}_{v,w}=0$ for $v\neq w$. Since $\Omega^{(2)}$ is symmetric, without loss of generality let $w$ be a non-descendant of $v$. Then, so are all of its parents. Therefore,
    \begin{align*}
        \Omega^\inParens{2}_{v,w}
        &= \inBrackets{\inParens{S - \Lambda^T S} - \inParens{S - \Lambda^T S}\Lambda}_{v,w}
        \\
        &= \inParens{S_{v,w} - \sum_{u\in\pa(v)} S_{u,w} \lambda_{u,v}}
        - \sum_{y\in\pa(w)} \inParens{S_{v,y} - \sum_{u\in\pa(w)} S_{u,y} \lambda_{u,v}} \lambda_{y,w}
        \\
        &= 0.
    \end{align*}

    The same observation applies to $\Omega^{(3)}$ defined in equation~(\ref{equ diagonalization of T}) as follows. 
    Let $z \in V$, and $w$ be a non-descendant of $v$, in which case, so are all of its parents $y\in\pa(w)$. 
    By symmetry of $\Omega^\inParens{3}$, 
    we need only compute $\Omega^\inParens{3}_{v,w,z}$,
    \begin{align*}
        \Omega^\inParens{3}_{v,w,z}
        =& \inBrackets{\inParens{\inParens{T - T\tucker\Lambda} - \inParens{T - T\tucker\Lambda}\tucker\Lambda}\tucker(I-\Lambda)}_{v,w,z}
        \\
        =& \sum_{x\in V}\left( \inParens{T_{v,w,x} - \sum_{u\in\pa(v)} T_{u,w,x} \lambda_{u,v}} \right.
        \\
        &\left.- \sum_{y\in\pa(w)} \inParens{T_{v,y,x} - \sum_{u\in\pa(w)} T_{u,y,x} \lambda_{u,v}} \lambda_{y,w} \right) \inBrackets{I-\Lambda}_{x,z}
        \\
        =& \sum_{x\in V}\left( 0 - \sum_{y\in\pa(w)} 0 \cdot \lambda_{y,w} \right) \inBrackets{I-\Lambda}_{x,z}
        \\
        =& 0.
    \end{align*}
    Thus, every entry but the diagonals of $\Omega^\inParens{2}$ and $\Omega^\inParens{3}$ is $0$.
\end{proof}

\begin{remark}
    In Lemma~\ref{lem matrices recover the graph}, we require that each $M_v$ defined in equation~\eqref{equ definition of M_v} drops rank.
    However, we still get the result of Lemma~\ref{lem matrices recover the graph} by instead ensuring that the smaller matrices,
    \begin{equation*}
        M_v' := \begin{bmatrix}
                S_{\pa(v),\nd(v)} &T_{\pa(v),N_v} 
                \\
                S_{v,\nd(v)} &T_{v,N_v}
            \end{bmatrix}
    \end{equation*}
    drop rank.
    Here the set $N_v \subseteq \flatten{\nd(v)}{V}$ is defined as follows,
    \begin{equation*}
        N_v := \{(w,z) \vert 
                w \in \nd(v) \textnormal{ and }
                z \in \nd(w) \cup \{v,w\}
            \}.
    \end{equation*}
\end{remark}
\begin{proof}
    This follows from the step in the poof at which we claimed that 
    "by symmetry of $\Omega^\inParens{3}$, we need only compute $\Omega^\inParens{3}_{v,w,z}$."
    A more detailed analysis shows that we can restrict further to $(w,z) \in N_v$ as claimed.
\end{proof}

\subsection{Proof of Proposition~\ref{prop other matrices}}

    Let $v \in V$ and $A \subseteq V\setminus v$.
    Recall that for any $w \in \nd(A\cup v) \subseteq \nd(v)$ and $z \in V$,
    \begin{align}
        s_{vw} 
        =& \sum_{u\in\pa(v)} \lambda_{uv} s_{uw},
        \label{equ:moreEquations_s}
        \\
        t_{vwz} 
        =& \sum_{u\in\pa(v)} \lambda_{uv} t_{uwz}.
        \label{equ:moreEquations_t}
    \end{align}
    If $u \in \pa(v)\cap A$, then $w\in\nd(u)$, so we have
    \begin{align}
        s_{uw} 
        =& \sum_{x\in\pa(u)} \lambda_{xu} s_{xw},
        \label{equ:moreEquations_s_parent}
        \\
        t_{uwz} 
        =& \sum_{x\in\pa(u)} \lambda_{xu} t_{xwz}.
        \label{equ:moreEquations_t_parent}
    \end{align}
    We see that we can replace all terms indexed by $u \in \pa(v)\cap A$ in equations~\eqref{equ:moreEquations_s} and~\eqref{equ:moreEquations_t} using equations~\eqref{equ:moreEquations_s_parent} and~\eqref{equ:moreEquations_t_parent}.
   
    \begin{align*}
        s_{vw} 
        =& \sum_{x\in P_{1,v}} \lambda_{xv}^\inParens{1} s_{xw},
        \\
        t_{vwz} 
        =& \sum_{x\in P_{1,v}} \lambda_{xv}^\inParens{1} t_{xwz},
    \end{align*}
   where $\lambda_{xv}^\inParens{1}$ is some polynomial in the $\lambda_{ij}$'s, and the set $P_{1,v} \subseteq V$ is defined as follows.
    For each $u\in\pa(v)$, if $u \notin A$, then include $u$ in $P_{1,v}$. Otherwise if $u \in A$, then include $\pa(u)$ in $P_{1,v}$.

   Repeating this procedure, for each $n\geq 2$, we obtain sets $P_{n,v}$ and equations
    \begin{align*}
        s_{vw} 
        =& \sum_{x\in P_{n,v}} \lambda_{xv}^\inParens{n} s_{xw},
        \\
        t_{vwz} 
        =& \sum_{x\in P_{n,v}} \lambda_{xv}^\inParens{n} t_{xwz},
    \end{align*}
    where the set $P_{n,v}$ is defined as follows. For each $u \in P_{n-1,v}$, if $u \notin A$, then $u \in P_{n,v}$, else if $u \in A$ then $P_{n,v}$ contains $\pa(u)$.

    Since $G$ is acyclic, there is no infinite chain of vertices $(u_n)$ such that $u_{n+1}\in\pa(u_n)$,
    thus the sequence of sets $P_{n,v}$ terminates.
    The terminating set contains only elements of $V$ which are parents of either $v$ or $A$ and do not lie in $A$.
    Thus, is contained in $\pa(\{v\}\cup A)$.
    Hence, the equations have the form
    \begin{align*}
        s_{vw} 
        =& \sum_{x\in\pa(\{v\}\cup A)} \lambda_{xv}^\inParens{\infty} s_{xw},
        \\
        t_{vwz} 
        =& \sum_{x\in\pa(\{v\}\cup A)} \lambda_{xv}^\inParens{\infty} t_{xwz},
    \end{align*}
    for some coefficients $\lambda_{xv}^\inParens{\infty}$, some of which may be $0$.
   This implies that the last row of
    \begin{equation*}
        \begin{bmatrix}
                R_{\pa(\{v\}\cup A),\nd(\{v\}\cup A)}
                \\
                R_{v,\nd(\{v\}\cup A)}
        \end{bmatrix},
    \end{equation*}
    is a linear combination of its other rows.

\section{Applications}\label{sec:application}
In this section we apply our results to the problem of recovering a sink node, and we examine the threshold sensitivity of the condition numbers which we use for determining matrix ranks.

We sample from a LiNGAM with graph $G = (V,E)$ and ask the question: "Can we recover a sink node from the data?".
We first compute estimators $\hat S$ and $\hat T$ for $S$ and $T$, respectively.
For each vertex $i \in V$, we create the following matrix,
\begin{equation}
    \hat M_i = 
    \begin{bmatrix}
        \hat S_{(V\setminus i), (V\setminus i)} &\hat T_{(V\setminus i), (V\setminus i)\times(V\setminus i)}
        \\ \hat S_{i,(V\setminus i)} &\hat T_{i,(V\setminus i)\times(V\setminus i)}
    \end{bmatrix}.
    \label{equ:matrix_estimator}
\end{equation}
By Theorem~\ref{thm:main}, the true matrix $M_i$ will drop rank if and only if $i$ is not an ancestor of any other vertices in the graph.
Then for each $i$ we compute the singular value decomposition (SVD) of $M_i$ and return its condition number $c_i$.
We obtain a guess for the sink of our unknown graph $G$ by picking the vertex $\hat i = i$ for which the singular value $c_i$ is minimal.
That is, for which the matrix $M_i$ is most likely to drop rank.


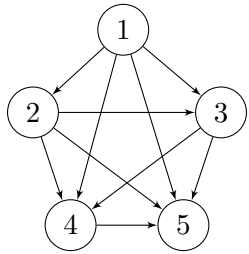
\begin{figure}
    \centering
    \begin{tikzpicture}
        \tikzset{vertex/.style = {shape=circle,draw,minimum size=1.5em}}
        \tikzset{edge/.style = {->,> = latex'}}
        \node[vertex] (1) at  (1.7, 1.1) {$1$};
        \node[vertex] (2) at  (0.5,   0) {$2$};
        \node[vertex] (3) at  (  3,   0) {$3$};
        \node[vertex] (4) at  (  1,-1.5) {$4$};
        \node[vertex] (5) at  (2.5,-1.5) {$5$};
        \draw[edge] (1) to (2);
        \draw[edge] (1) to (3);
        \draw[edge] (1) to (4);
        \draw[edge] (1) to (5);
        
        \draw[edge] (2) to (3);
        \draw[edge] (2) to (4);
        \draw[edge] (2) to (5);

        \draw[edge] (3) to (4);
        \draw[edge] (3) to (5);

        \draw[edge] (4) to (5);
    \end{tikzpicture}
    \caption{
        The complete DAG on $5$ vertices
    }
    \label{fig 5 node acyclic}
\end{figure}
\begin{figure}
    \centering
    \includegraphics[width=0.75\linewidth]{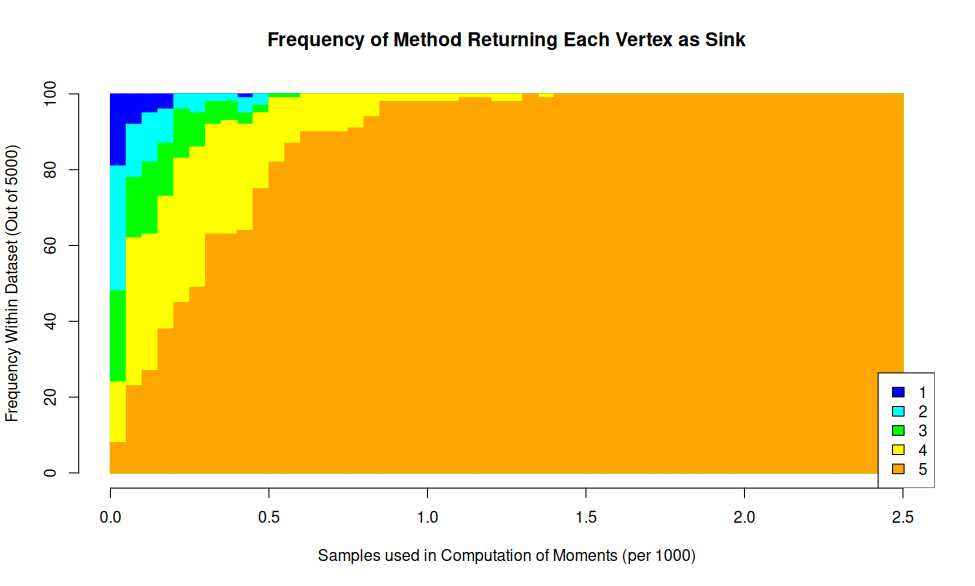}
    \caption{
        At each number of samples $n$, 
        we record the number of times each vertex was assigned as the sink out of the entire data set.
        The graph is a bar chart, each bar represents $200$ data points.
    }
    \label{fig:Graph}
\end{figure}

Let $G$ be the graph given in Figure~\ref{fig 5 node acyclic}.
To produce Figure~\ref{fig:Graph}, we perform $50$ runs, where each error is sampled from a $\Gamma(5,1)$ distribution. 
Using $100$ different sample sizes, we compute $\hat S$ and $\hat T$,
and produce an index $\hat i \in \{1,\dots 5\}$, an estimate of a sink node. For each vertex, we record the number of times it was output and the number of samples used in the computation.

The result is a graph, Figure~\ref{fig:Graph}, which at each number of samples $n$, 
records the proportion of trials on $n$ samples in which the method returned each vertex as the sink. 

As the number of samples grows, we observe that the vast majority of times we obtain the correct sink node (node 5). A small proportion of the experiments wrongly label node 4 as a sink node, most likely due to numerical error.

\subsection*{Sensitivity to Changes of Threshold Values}

\begin{figure}
    \centering
    \begin{tikzpicture}
        \tikzset{vertex/.style = {shape=circle,draw,minimum size=1.5em}}
        \tikzset{edge/.style = {->,> = latex'}}
        \node[vertex] (1) at  (  0, 1.5) {$1$};
        \node[vertex] (2) at  (  0,   0) {$2$};
        \node[vertex] (3) at  (1.5,   0) {$3$};
        \draw[edge] (1) to (2);
        \draw[edge] (2) to (3);
    \end{tikzpicture}
    \caption{
        The line DAG on $3$ vertices
    }
    \label{fig 3 node line}
\end{figure}
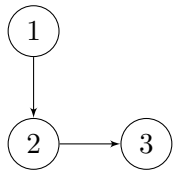
\begin{figure}
    \centering
    \includegraphics[width=0.5\linewidth]{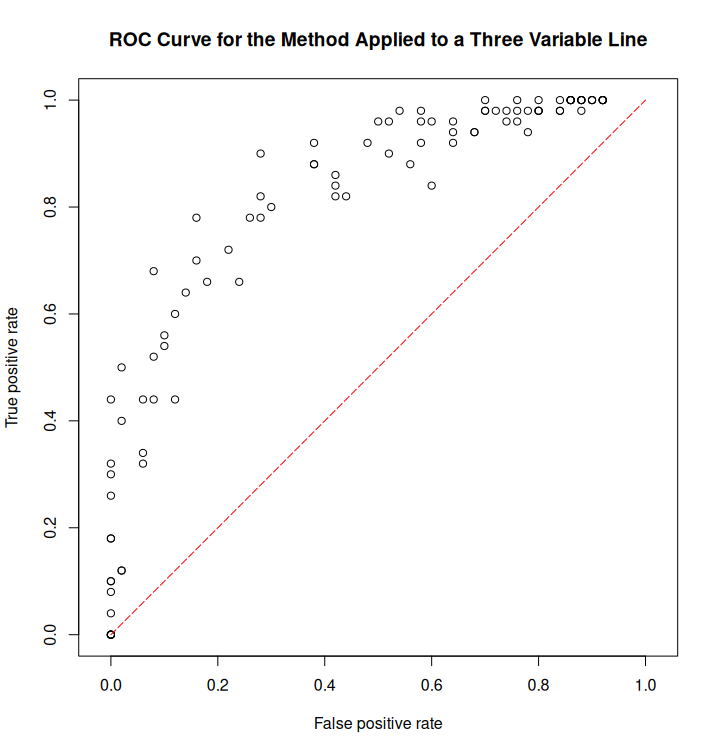}
    \caption{
        True positive rate plotted against the false positive rate for different threshold values.
        We use $5000$ samples to construct $\hat S$ and $\hat T$ for each run, 
        Each data point represents one threshold value and uses $50$ runs to approximate the true and false positive rates.
    }
    \label{fig:roccurve}
\end{figure}
As above, we sample from a LiNGAM with graph given in Figure~\ref{fig 3 node line}.
Picking a threshold value $t$, 
we guess that a vertex $v$ is a sink of $G$ if the rank of the matrix on line~\ref{equ:matrix_estimator} is less than $t$.
The rate at which this method correctly chooses $3$ as a sink
versus choosing either of $1$ or $2$ as a sink for different values of $t$ is given in Figure~\ref{fig:roccurve}.
The number of samples used to estimate $S$ and $T$ is fixed constant at $5000$. Note that the method could label more than one vertex as a sink.

Figure~\ref{fig:roccurve} shows a nicely shaped ROC curve, which would allow us to pick an appropriate threshold value.
For instance, if we wanted to ensure a false positive rate of at most $20\%$, then we can pick a threshold with which we would expect a true positive rate between $60$ and $80\%$.

\section{Discussion}~\label{sec:discussion}

In this paper we studied the set $\mathcal M^{\leq 3}(G)$ of second- and third-order moments $S$ and $T$ of a random  vector which satisfies  a linear structural  equation model with respect to a directed acyclic graph $G$. We derived explicit polynomial constraints in the entries of $S$ and $T$ which cut out the set $\mathcal M^{\leq 3}(G)$ (Theorem~\ref{thm:main}). We then  proved  that our equations generalize earlier work in  the case of polytrees (Section~\ref{sec:additional_equations}). Furthermore, we noted that when restricted to the covariance matrix $S$, our equations imply the Local Markov Property.

This work opens up more interesting questions in the field of algebraic statistics. The existing cumulant-based algorithm for recovering a latent variable LiNGAM of~\cite{schkoda2024causal} suggests that one could obtain a similar characterization of the set of second-, third-, and potentially higher-order cumulants in a latent variable LiNGAM. Furthermore,  even when all variables are observed if the underlying distribution of the coordinates of the random vector $X$ is symmetric around 0, using higher-order moments would be needed in order to uniquely recover the graph. We believe that extending our results to cumulants of order higher than 3 should be completely analogous to the proof of Theorem~\ref{thm:main}. 

While potentially more difficult, it would be quite interesting to study the model of second and third-order moments when the graph $G$ is allowed to have directed cycles. Finding the defining equations in this case would be quite useful in designing a causal discovery algorithm, extending previous work which only applies to cycle-disjoint graphs~\cite{drton2025cyclic}.

We also believe that the determinantal equations arising from Theorem~\ref{thm:main} could potentially shed light in the context of linear Gaussian hidden variable models which are not cut out by determinantal constraints, such as Verma constraints, the Pentad, and others~\cite{drton2020nested}. Adding the third-order moments in such models (and then eliminating them from the defining ideal) should shed light on how to obtain their characterization in general.

\section{Acknowledgments}
We would like to thank Mathias Drton for helpful discussions. Elina Robeva was supported by a Canada CIFAR AI Chair and an NSERC Discovery Grant (DGECR-2020-00338).

\bibliography{biblio}

\end{document}